\documentclass[10pt]{article}

\usepackage[utf8]{inputenc}
\usepackage[T1]{fontenc}
\usepackage[english, french]{babel}
\usepackage{amsmath, amsfonts, amstext, amssymb} 
 \usepackage{color,colortbl}
 \usepackage{fancyhdr}
 \usepackage{xcolor,pict2e}
 \definecolor{myaqua}{rgb}{0.0,0.5,0.55}
 \definecolor{lightaqua}{rgb}{0.75,0.95,0.95}
\usepackage{ epsfig, lscape}
 \usepackage[colorlinks = true,
            linkcolor = myaqua,
            urlcolor  = blue,
            citecolor = red,
            pdfpagemode=UseOutlines,
            bookmarksnumbered=true,
            pdfpagelabels,
            breaklinks]{hyperref}

\usepackage{tikz}
\usetikzlibrary{shapes}
\usepackage{lipsum}

\usepackage{rcs}

\usepackage{graphicx}
\newcommand{\N}{\mathbb{N}}
\newcommand{\Z}{\mathbb{Z}}
\newcommand{\R}{\mathbb{R}}
\newcommand{\M}{\mathcal{M}}
\newcommand{\A}{\mathcal{A}}

\newcommand{\F}{\mathfrak{P}}
\newcommand{\C}{\mathfrak{C}}

\newtheorem{theo}{Theorem}

\newtheorem{lem}{Lemma}
\newtheorem{prop}{Proposition}

\newtheorem{rem}{Remark}

\newenvironment{preu}{\textbf{Proof }\rm}{\par\hfill$\square$}
\newenvironment{preu 1}{\textbf{Proof of Theorem \ref{t2}.}\rm}{\par\hfill$\square$}
\newenvironment{preu 2}{\textbf{Proof of Theorem \ref{t1}.}\rm}{\par\hfill$\square$}
\newenvironment{preu 3}{\textbf{Proof of Theorem \ref{th3}}\rm}{\par\hfill$\square$}
\newenvironment{preu 4}{\textbf{Proof of Theorem \ref{theo3}}\rm}{\par\hfill$\square$}
\newenvironment{preu 5}{\textbf{Proof of Lemma \ref{lmm3} }\rm}{\par\hfill$\square$}
\newtheorem{defi}{Definition}
\newtheorem{exple}{Example}

\usepackage{datetime}
\usepackage{lineno}
\usepackage{floatrow}

\def\lin#1#2{\textcolor[rgb]{0.6,0.6,0.6}{\vspace*{#1mm} \hrule
   height 3 pt \vspace*{#2mm}}}
\def\bt{\begin{tabular}}
\def\et{\end{tabular}}
\def\and{\mbox{ and }}

\def\1{{\bf 1}}

 \def\boxx#1#2#3#4#5{
 {\linethickness{#4pt}\put(#1,#5){\color{myaqua}{\line(1,0){#3}}}}
 \multiput(#1,#2)(0,#4){2}{\line(1,0){#3}}
 \multiput(#1,#2)(#3,0){2}{\line(0,1){#4}}
  }

\begin{document}


 $\mbox{ }$

 \vskip 12mm
 
 { 

{\noindent{\Large\bf\color{myaqua}
   Contribution on the Intrinsic Ergodicity of the Negative Beta-shift  }} 
%
\\[6mm]
{\bf Florent NGUEMA NDONG }}
\\[2mm]
{ 
 $^1$ Universit\'e des Sciences et Techniques de Masuku 
 \\
Email: \href{mailto:florentnn@yahoo.fr}{\color{blue}{\underline{\smash{florentnn@yahoo.fr}}}}\\[1mm]

\lin{5}{7}

 {  
 {\noindent{\large\bf\color{myaqua} Abstract}{\bf \\[3mm]
 \textup{
  Let $ \beta $ be a real number less than -1. In this paper, we prove the uniqueness of the measure with maximal entropy of the negative $\beta$-shift. Endowed with the shift, this symbolic dynamical system is coded under certain conditions, but in all cases, it is shown that the measure with maximal entropy is carried by a support coded by a recurrent positive code. One of the difference between the positive and the negative $\beta$-shift is the existence of gaps in the system for certain negative values of $ \beta $. These are intervals of negative $\beta$-representations (cylinders) negligible with respect to the measure with maximal entropy, which is a measure of Champernown.   
 }}}
 \\[4mm]
 {\noindent{\large\bf\color{myaqua} Keywords}{\bf \\[3mm]
  Negative basis; $\beta$-expansions; coded system; transitivity; intrinsic ergodicity 
}}
\\[4mm]
 {\noindent{\large\bf\color{myaqua} 2010 MSC}{\bf \\[3mm] 11K16; 11B05; 37A05; 37A25; 37B10
}}}
\lin{3}{1}

\renewcommand{\headrulewidth}{0.5pt}
\renewcommand{\footrulewidth}{0pt}

 \pagestyle{fancy}
 \fancyfoot{}
 \fancyhead{} 
 \fancyhf{}
 \fancyhead[RO]{\leavevmode \put(-105,0){\color{myaqua}F. NGUEMA NDONG} \boxx{15}{-10}{10}{50}{15} }
 \fancyfoot[C]{\leavevmode
 \put(-2.5,-3){\color{myaqua}\thepage}}

 \renewcommand{\headrule}{\hbox to\headwidth{\color{myaqua}\leaders\hrule height \headrulewidth\hfill}}

 \section{Introduction}

{ 
\selectfont
 \noindent

Consider a real number $ \beta $ with modulus greater than 1. Since the seminal paper \cite{MR0097374} on expansions of numbers in non-integer positive base, many mathematicians have become interested in the properties of the $\beta$-shift with $ \beta > 1 $. For instance, in \cite{MR939059}, the author established the intrinsic ergodicity. The condition to have the specification property is given in \cite{MR545668}. The concept of intrinsic ergodicity was first explored by W. Parry in \cite{MR0142719}: it is about the uniqueness of a measure with maximal entropy, and the study of the indecomposability of the system into several invariant subsets non negligible with respect to this measure. Note that the set of invariant measures on a dynamical system is non-empty.

Vittorio Grünwald discovered the negabinary system. It is a non-standard positional numeral system with the unusual property that negative and positive numbers can be represented without a sign bit. It was been used in the experimental Polish computers SKRZAT 1 (see for example \cite{FIETT19601250}) and BINEG in 1950. In 2009 S. Ito and T. Sadahiro, in the seminal paper \cite{MR2534912} extended the notion of representation of numbers without a sign to all negative bases. One of interesting features of the theory of expansions of numbers is the link it creates between analytic number theory and symbolic dynamics.  
In the last decade, various papers had been devoted to the similarities and differences between the positive and negative $ \beta$-transformations, and basic properties of the strings corresponding to negative $\beta$-expansions had been derived. This paper fits in this line of research. More precisely, we study the intrinsic ergodicity of the negative $\beta$-shift. The main result is stated in Theorem \ref{thIE}. The question on the intrinsic ergodicity of the negative $\beta$-shift was first explored by S. Mao and Y. Kenichiro in \cite{article}. The authors proved that if $ \beta < -\frac{1+\sqrt{5}}{2} $, the $ \beta$-shift is intrinsically ergodic. Theorem \ref{thIE} of this paper shows that this result can be extended to all values $ \beta $ less than -1. For $ \beta $ taken in $ (-\frac{1+\sqrt{5}}{2}, -1 )$, the negative $ \beta$-shift (denoted by $S_{\beta}$) is not transitive. Thus, in this case, a measure with maxial entropy is carried by a subset strictly included in the system. In fact, the support of a measure with maximal entropy is transitive.
We exhibit this support by distinguishing each transitive sub-shift contained in $ S_{\beta} $ with an associated code (prefix or suffix).
Our approach relies first and foremost on well-known result on coded systems of works of G. Hansel and F. Blanchard in \cite{MR858689} and A. Bertrand (see \cite{MR939059}). To prove the intrinsic ergodicity for all values of $\beta$, we proceed in two steps. Firstly, the measures with maximal entropy of the negative $\beta$-shift are carried by the same support. Secondly, this support, endowed with the shift, has a unique measure with maximal entropy.

Endowed with the shift, the positive and negative $\beta$-shifts share many properties. On the other hand, they differ in many aspects. For instance, the positive $\beta$-shift is coded for all values of $ \beta$ ($\beta>1$). In contrast, the negative $\beta$-shift is coded if and only if $ \beta$ is less than or equal to $ -\frac{1+\sqrt{5}}{2} $ and the $ \beta$-expansion of the left end-point of the domain of the negative $\beta$-transformation is not periodic with odd period. Nevertheless, the measure with maximal entropy of the negative $\beta$-shift is a Champernown measure of a recurrent prefix (or suffix) code and it is mixing (see Theorem \ref{pp1} and Theorem \ref{Mix}). Moreover, a strange phenomenon labels the difference between positive and negative $\beta$-shifts: the existence of gaps in the system for the negative case. It is about subsets of the negative $\beta$-shift negligible with respect to the measure with maximal entropy. We exhibit these intervals of sequences (cylinders). Such cylinders are carried by intransitive words (see  Theorems \ref{InW1} and \ref{InW2}). In the domain of the negative $\beta$-transformation (which we will denote by $I_{\beta}$), this phenomenon was investigated in \cite{MR2974214}. 
 
  The contents of this paper are as follows. We start our study by generalities on symbolic dynamical system. More precisely, we begin by a brief overview of coded systems, the notion of tower of a prefix code (introduced by G. Hansel and F. Blanchard in \cite{MR858689}) and $\beta$-shift. 
 The second part of the paper is devoted to the intrinsic ergodicity. We start this section by recalling the codes of the possible supports of a measure with maximal entropy of the negative $\beta$-shifts. Next, we give intransitive words. At last, we determine the measure with maximal entropy and prove its uniqueness.

\section{Generality}

In this section, we briefly recall several facts about coded systems and the representation of numbers in real bases. 

\subsection{Coded System}

Let $\A$ be an alphabet, $(X,T)$ a symbolic dynamical system on $ \A$ and $ L_X$ the associated language. In the following, we denote by $ \A^*$ the free monoid generated by $ \A $ and $ \A^{+} = \A^{*}\backslash\{ \varepsilon \} $ where $ \varepsilon $ is the empty word, $ \N^{\times} = \N\backslash\{0\} $, $\R^{\times} = \R\backslash\{0\}$ and $ \Z^{\times} = \Z\backslash\{0\}$. We recall some definitions given in \cite{MR858689}.

\begin{defi}
A language $L$ is said to be transitive if for all pair of words $(u,v)$ of $L^2$, there exists $w$ in $ \mathcal{A}^{*} $ such that  $uwv$ belongs to $L$.
\end{defi}

A symbolic dynamical system will be said transitive if the associated language is transitive. More generally,consiering a topological dynamical system $ (X,T)$, if for all open sets $U$ and $ V$ of $ X $, one has $ U \cap T^{-n}V \neq \varnothing $ for some $ n$ in $\Z$. This is equivalent to saying that the orbit $ \bigcup \limits_{n \in \Z} T^n U $ of all non-empty open set $ U $ of $X$ is dense in $ X $ (see \cite{MR0352411}).

\begin{defi}
A code $Y$ on $\A$ is a language such that, for any equality 
\begin{equation}
x_1x_2\cdots x_n = y_1y_2 \cdots y_k,
\end{equation}
for any $ x_i$, $y_j \in Y$, one has $ x_i=y_i $ and $ k = n$. 
\end{defi}

A prefix (resp. suffix) code is a language $ \C $ of $\A^{+}$ for which no word is the prefix (resp. suffix) of another. That is,
\begin{equation}
 \forall u, v \in \C, u = vw \Rightarrow u = v \text{ and } w = \varepsilon.
\end{equation}

\begin{defi}
A coded system is a sub-shift for which there exists a (non-unique) language $ Y $ such that $ S $ is the closure of $ Y^{\infty}$.
\end{defi}

Let $ L $ be a language on an alphabet $ \A $. The radius $ \rho_L $ of the power series $ \sum\limits_{n \geq 1} card (L\cap \A^n)z^n $ is called the radius of convergence of $ L $. 

\begin{defi}
Let $ X $ be a symbolic dynamical system and $ a_1 a_2 \cdots a_k \in L_X$. We denote by $ _{m}[a_1\cdots a_k ]$, the set of sequences $ (x_i)_{i \in \Z}$ of $ X $ such that 
\begin{equation}
x_m x_{m+1} \cdots x_{m+k-1} = a_1 a_2 \cdots a_k.
\end{equation}
This set is called cylinder carried by $ a_1a_2\cdots a_k $ at $ m $.

For $ m = 0 $, in order to simplify the notations, we will denote  $_{0}[x] $ by $ [x] $.
\end{defi}

Consider a symbolic dynamical system $ X $ on an alphabet $ \A $ and $ x \in L_X $. The length of a word $ x $, denoted by $ l(x) $, is the number of letters of $ x $. 
\begin{equation}
[x = x_1 x_2 \cdots x_k, \text{ with $ x_i \in \A $ }]\Rightarrow l(x) = k.
\end{equation}

\begin{defi}
A prefix code $ \C $ is said to be recurrent positive if 
\begin{equation}
\sum\limits_{x \in \C} \rho_{\C^{*}}^{ l( x ) } = 1 \text{ and } \sum\limits_{ x \in \C} l( x)  \rho_{\C^{*}}^{l( x )}<+\infty.
\end{equation}

If we set $ \rho_{\C^{*}} = \frac{1}{\beta}$, one has $ 1 = \sum\limits_{n \geq 1} \frac{c_n}{\beta^n} $ and $ \sum\limits_{n \geq 1}\frac{nc_n}{\beta^n} < +\infty $ where $ c_n $ counts the number of words of length $ n $ in $ \C $ and $ \C^{*} $ the free monoid generated by $ \C $.
\end{defi} 

\begin{rem}\label{Nrem1}
From Proposition 2.15 of \cite{MR858689}, if $\mu$ is an invariant measure on $\C^{\Z}$ endowed with the shift and $h(\mu)$ is its entropy, then, 
 \begin{itemize}
  \item[(1)] one has 
  \begin{equation}
   h(\mu) \leq - l(\C, \mu)\log \rho_{\C^{*}}
   \label{(1)}
  \end{equation}
  where
  \begin{equation}
  l(\C, \mu) = \sum\limits_{ x \in \C }l(x)  \mu([x]).
  \end{equation}
\item[(2)] In \eqref{(1)}, the equality holds if one has both following conditions:
  \item[(a)] $ \sum \limits_{x \in \C} \rho^{l( x )}_{\C^{*}} = 1 $,
  \item[(b)] $ \mu $ is a probability of Bernoulli on $\C^{*}$ defined by
  $ \mu([x])=\rho^{l(x)}_{\C^{*}}$, $ x \in \C $.
 \end{itemize}
\end{rem}
 
\begin{defi}
A measurable topological dynamical system $ (X, m, g) $ is said to be \textit{ergodic} if for any measurable $g$-invariant set $ B \subset X $, one has $m(B) = 0$ or $ m(B) = 1$. One also says that $ m $ is an ergodic measure with respect to $ g $ or $ g $ is ergodic with respect to $ m $. It is said to be mixing if for all measurable sets $ A $ and $ B $,
 \begin{equation}
  \lim \limits_{n \rightarrow +\infty} m(g^{-n}(A) \cap B) = m(A)m(B).
 \end{equation}
\end{defi}

\subsection{Tower associated to a prefix code}

More details on the notion of the tower of a prefix code can be found in \cite{MR858689}.
Let $ \Omega $ be the subset of $ \C^{\Z}\times \N $ such that :
\begin{equation}
 ((x_n)_{n \in \Z}, i) \in \Omega \Rightarrow 1 \leq i \leq l( x_0 ).
\end{equation}
We can identify $ (x_n)_{n \in \Z} $ with an element $ x $ of $ \C^{\Z}$ which is a concatenation of words $ x_i $ of the code $\C$.

We define a map $ T $ from $ \Omega $ into itself by:
\begin{equation}
 T((x_n)_{n \in \Z}, i) = \begin{cases}
                           ((x_n)_{n \in \Z}, i+1)   & \text{ if $ i < l( x_0 ) $ } \\
                           ((x_{n+1})_{n \in \Z}, 1) & \text{ if $ i = l( x_0 ). $ }
                          \end{cases}
\end{equation}
The pair $ (\Omega, T) $ is called \textit{the tower associated} to $ \C $. We have the two important following facts (see \cite{MR858689}).
 \begin{itemize}
  \item[•] When $ \overline{\mu} $ ranges through $\M_T(\Omega)$, $\sup_{\overline{\mu}}h(\overline{\mu}) = -\log \rho_{\C^{*}}$. 
  
  \item[•] There exists one and only one invariant probability $ \mu $ on $ \Omega $ such that $ h(\overline{\mu})=-\log \rho_{\C^{*}} $ if only if $\C$ is recurrent positive. In this case, $ \overline{\mu} $ is the unique invariant probability (and thus ergodic) on $ \Omega $ inducing on $ \C^{\Z} $ a probability $\mu$ of Bernoulli defined by:
   \begin{equation}
   \mu ( [ x ] )= \rho^{l(x)}_{\C^{*}}, \text{ $ x \in \C $ }.
  \end{equation}
 \end{itemize}
The first statement is usually referred to as the variational principle. When the dynamical system is coded by a recurrent prefix code $\C$, the set of periodic points is dense. One says that two words $ u $ and $ v $ of the language $L$ are in the same class of syntactic monoid if for all pair of words $(a,b)$,
\begin{equation}
 aub \in L \Longleftrightarrow a v b \in L. 
\end{equation}
A symbolic dynamical system $S$ is said to be rational (or sofic) if the number of classes of language associated to is finite. In this case, it is coded by a recurrent positive prefix code.
The induced measure $ \overset{\bullet}{\mu} $ on the coded system by the measure on the tower is particularly simple: if $ u_1u_2 \cdots u_r $ is a word of $ L $,
\begin{equation}
 \overset{\bullet}{\mu}([u_1u_2 \cdots u_2]) = \sum\limits_{ m } \mu([m]) = \dfrac{1}{\sum\limits_{n \geq 1}\frac{n.c_n}{\beta^n}}
 \sum \limits_{ m } \frac{1}{\beta^{l( m )}}
\end{equation}
where in the sums, $ m $ is all words of the form $ m = a u_1u_2 \cdots u_k b $; $ a $ is a proper prefix of a word of the code and $ b $ is a proper suffix or the empty word.

Each class of the syntactic monoid is associated to a positive constant $ \lambda $ such that if $ u $ belongs to this class, the cylinder $[u]$ has a measure
\begin{equation}
 \overset{\bullet}{\mu}([u]) = \frac{\lambda}{\beta^{l(u )}}.
\end{equation}
So, the measure of a cylinder depends on the class of its support (the word $ u $) and its length. If there exist two positive numbers $ \epsilon $ and $ M $, with  
 $ 0 < \epsilon < \lambda < M $ for all classes, the measure is said to be homogeneous. The ratio between measures of cylinders of the same length is controlled by $ \frac{\epsilon}{M} $ and $ \frac{M}{\epsilon} $. If the greatest common divisor (G.C.D) of lengths of words of the code is 1, then the measure is mixing. In this case, it is said to be of Champernowne.

\subsection{Beta-shift}

Consider a real number $ \beta $ with modulus greater than 1. Let us approach the question of representing of numbers using powers of $ \beta $ from a more general point of view. Let 
\begin{equation}
l_{\beta} = \begin{cases}
          0 &\text{ if $ \beta > 1 $ } \\
          \frac{\beta}{1-\beta} &\text{ if $ \beta < -1 $ }
           \end{cases} \text{ and $ r_{\beta} = l_{\beta} + 1 $. }
\end{equation} 

Define the map $ T_{\beta}$ from $I_{\beta} = [l_{\beta}, r_{\beta} )$ into  itself by  
\begin{equation}
  T_{\beta}(x) = \beta x - \lfloor \beta x -l_{\beta} \rfloor.
\end{equation}
The expansion in base $ \beta $ of a real $ x $ (denoted by $ d(x, \beta )$) is given by the following algorithm.
We find the smallest non negative integer $ n $ for which one has $ \frac{x}{\beta^n} \in I_{\beta} $. The $\beta$-expansion is given by the sequence $ d(x, \beta) =x_{-n+1}\cdots x_0 \cdot x_1x_2 \cdots $ such that
\begin{equation}
  x_{-n+i} = \lfloor \beta T_{\beta}^{i-1}\left(\frac{x}{\beta^n}\right)-l_{\beta} \rfloor,\text{ $ i \geq 1 $ }.
\end{equation}
Let $ (r_i^*)_{i \geq 1 } $ be the sequence of digits defined by:
\begin{equation}
r^*_1 = \lfloor \beta r_{\beta} - l_{\beta} \rfloor
\end{equation}
and for all integers $ i \geq 2 $,
\begin{equation}
\begin{aligned}
r^*_i &= \lfloor \beta T_{\beta}^{i-2}(\beta r_{\beta}-r^*_1)-l_{\beta} \rfloor\\
    &=\lfloor \beta^i r_{\beta}- \sum\limits_{k=1}^{i-1} r^*_k \beta^{i-k}-l_{\beta} \rfloor.
\end{aligned}
\end{equation}
Consider an alphabet $ \A $ endowed with an order $ \prec_{\delta}$ such that for all sequences of digits $(x_i)_{i \geq 1} $ and $(y_i)_{i \geq 1}$ over $ \A $, 
\begin{equation}
(x_i)_{i \geq 1} \prec_{\delta} (y_i)_{i \geq 1} \Leftrightarrow \exists k \in \N^{\times}| \hspace{0.2 cm} x_i=y_i \hspace{0.2cm}\forall i < k, \hspace{0.5cm}\delta^k(x_k-y_k)< 0
\end{equation}
where $ \delta $ is the sign of $ \beta $.
\begin{itemize}
\item If $ \beta < -1 $, $ \prec_{\delta} $ is the alternating order (see for example \cite{MR2534912}, \cite{NguemaNdong20161}, \cite{NguemaNdong2019}).
\item If $ \beta> 1$, $ \prec_{\delta} $ is the classical lexicographic order on words. 
\end{itemize}

Let $ d(l_\beta, \beta) = \cdot d_1 d_2 \cdots $. 
Denote by $ f_{\beta} $ the map on words defined by: 
\begin{equation}
f_{\beta} ((x_i)_{i \in \Z}) = \sum\limits_{ k \in \Z} x_k \beta^{-k}.
\end{equation} 
If $d(x, \beta ) = (x_i)_{i \geq n }$, then 
\begin{equation}
x = f_{\beta} ((x_i)_{i \geq n})
\end{equation} 
And thus $ f_{\beta} ((r^{*}_i)_{i \geq 1}) = r_{\beta} $ and $ f_{\beta} ((d_i)_{i \geq 1}) = l_{\beta} $. 
Let $ (r_i)_{i \geq 1}$ be the sequence of digits such that
\begin{equation}
 (r_{i})_{i \geq 1} = \begin{cases}
                         \overline{(r_1^*, \cdots, r_{n-1}^*, r_{n}^*-1)} &\text{ if $ (r_i^*)_{i \geq 1} = (r_1^*, \cdots, r_{n}^*,d_1,d_2,\cdots ) $, $\beta>1$ }\\
                         \overline{(r_1^*, \cdots, r_{n-1}^*, r_{n}^*-1)} &\text{ if $ (r_i^*)_{i \geq 1} = (r_1^*, \cdots, r_{n}^*,d_1,d_2,\cdots ) $, $\beta<-1$ and $ n $ even }\\
                                  (r_i^*)_{i \geq 1} &\text{ otherwise } 
                      \end{cases}
\end{equation} 
and
\begin{equation}
 (d_{i}^*)_{i \geq 1} = \begin{cases}
                         \overline{(d_1, \cdots, d_{2n-2}, d_{2n-1}-1, 0)} &\text{ if $\beta< -1$, $ (d_i)_{i \geq 1} = \overline{(d_1, \cdots, d_{2n-1} )} $ }\\
                         (d_i)_{i \geq 1} &\text{ otherwise } 
                         \end{cases}\label{D}
\end{equation}      
where $ \overline{t} $ stands for infinite repetition of the string $ t $.
In fact, 
\begin{equation}
(d_i^*)_{i \geq 1} = \lim\limits_{x \rightarrow l_{\beta}^+}d(x, \beta) \text{ and } (r_i)_{i \geq 1} = \lim\limits_{x \rightarrow r_{\beta}^-}d(x, \beta).\label{lbrb}
\end{equation}
 Expansions in base $ \beta $ (or $ \beta$-expansions, $ |\beta|> 1 $) of real numbers are governed by the sequences $(d_i)_{i \geq 1}$ and $(r_i)_{i \geq 1}$.
A sequence $(x_i)_{i \geq 1} $ is the $ \beta $-expansion of a real $ x $ if and only if 
\begin{equation}
(d_i)_{i \geq 1} \preceq_{\delta}(x_{i+n})_{i \geq 1} \prec_{\delta}(r_i)_{i \geq 1}, \text{ $ \forall n \in \N $}.
\end{equation}
\begin{defi}
The $\beta$-shift $ S_{\beta}$ is the closure of the set of $\beta$-expansions.
\begin{equation}
S_{\beta} = \{ x_k x_{k+1} \cdots x_0 \cdot x_1\cdots \vert (d_i)_{i\geq 1} \preceq_{\delta}(x_i)_{i \geq m} \preceq_{\delta}(r_{i})_{i \geq 1}, \forall m\geq k\}.
\end{equation}
We also define the corrected $\beta$-shift $\tilde{S}_{\beta}$ as follows: 
\begin{equation}
\tilde{S}_{\beta} = \{ x_k x_{k+1} \cdots x_0 \cdot x_1\cdots \vert (d_i^*)_{i\geq 1} \preceq_{\delta}(x_i)_{i \geq m} \preceq_{\delta}(r_{i})_{i \geq 1}, \forall m\geq k, \forall k \}\label{sb}.
\end{equation}
\end{defi}
In fact, $ S_{\beta}$ and $ \tilde{S}_{\beta}$, endowed with the shift $\sigma$ have the same entropy. Moreover, each real has a representation in $ \tilde{S}_{\beta}$. Both bounds $(r_i)_{i \geq 1}$ and $ (d_i^*)_{i \geq 1}$ decide whether a digit string belongs to $ \tilde{S}_{\beta}$ or not. 
In the rest of the paper, $ L_{\beta} $ denotes the language of $ S_{\beta} $.

\begin{rem}\label{rmbta}
Let $ X $ be a symbolic dynamical system, $ \log t $ its entropy and $ H_n$ the number of words of $ L_X $ with length $ n $. Then, $ \frac{1}{t} $ is the smallest pole in modulus of $ \sum\limits_{ n \geq 0} H_n z^n $. 
\end{rem}

The intrinsic ergodicity of $S_{\beta}$ has been studied for $ \beta> 1$. Throughout the rest of the paper, we will be interested in the negative base case. In other words, $ \beta $ will be less than -1 and then, we will use the alternating order as the tool of controlability of words. In this case, the alternating order will be denoted by $ \preceq $ (or $\prec$) instead of $ \preceq_{\delta}$ (or $ \prec_{\delta}$).

\section{Intrinsic ergodicity of the negative beta-shift}

Let $ (X, T) $ be a topoloical dynaical system. The variational principle describes the relationship between topological entropy and Kolmogorov entropy of a measurable dynamical system. We denote by $ h_{top}(T) $ the topological entropy. If $ \mu $ is a $ T$-invariant measure of $ X $, we denote by $ h_{\mu}(T) $ the usual metric entropy of $ T $. In view of the variational principle, $ h_{top} $ coincides with the $ sup $ taken over the set of $ T$-invariant measures of the metric entropy. From \cite{MR0352411}, the $ sup $ can be considered just over the set of ergodic measures.
\begin{equation}
 h_{top}(T) = sup \{h_m(T) \vert m \in \M_T(X) \text{ is ergodic } \}.
\end{equation}
In symbolic dynamics, setting $ T= \sigma$ (the shift on words), the $sup$ exists because $\sigma$ is an expansive map (see more details in \cite{MR0457675}). That is we can find a real $ \theta > 0 $ such that for all $x$ and $y$ in any sub-shift $X$, $x\neq y$, there exists an integer $n $ for which one has:
\begin{equation}
 d( \sigma^{n}(x), \sigma^{n}(y)) \geq \theta 
\end{equation}
where $d$ is the metric defined by:
\begin{equation}
 d(x, y) = \sum \limits_{n \in \Z} 2^{-\vert n \vert} d(x_n , y_n) 
\end{equation}
 with $ x= (x_i)_{i \in \Z} $, $ y= (y_i)_{i \in \Z} $ and
\begin{equation}
 d(x_n, y_n)= \begin{cases}
               0 &\text{ if $ x_n = y_n $ } \\
               1 &\text{ if $ x_n \neq y_n $ }.
              \end{cases}
\end{equation}
It suffices to set $ \theta = \frac{1}{2} $.

Roughly speaking, a system is intrinsically ergodic if it has a unique measure (which is ergodic) of maximal entropy. The ergodicity of a measure implies the transitivity of its support. Thus, to study the intrinsic ergodicity of the negative beta-shift, we determine the possible supports of ergodic measures on $ S_{\beta}$ and then transitive subsystems of $ S_{\beta} $. It is well known that coded systems are transitive. Bearing this mind, and before adapting the study on the coded systems done in \cite{MR858689} and \cite{MR939059} to the negative beta-shift, we begin by proposing some codes we need to construct admissible words. 

Let $ (r_i)_{i \geq 1 } = \lim\limits_{x \rightarrow r_{\beta}}d(x, \beta) $. 
If $ \beta> 1 $, $(r_{i})_{i \geq 1}$ is the $ \beta$-expansion of 1. In this case, the $ \beta$-shift is coded by the language: 
\begin{equation*}
 \{ r_1 r_2 \cdots r_{n}j \hspace{0.1cm}| \hspace{0.1cm} 0 \leq j <r_{n+1 }, \text{ $ n \in \N $ }\}
\end{equation*}
with $ r_1\cdots r_n = \epsilon $ (the empty word) if $ n = 0 $.

Suppose $ \beta < -1 $ and let $ d(l_{\beta}, \beta) = (d_i)_{i \geq 1} $. In \cite{FloarXiv, NguemaNdong2019}, the author proved that $ S_{\beta}$ is coded if and only if  $ -\frac{1+\sqrt{5}}{2} \geq \beta $ and the $\beta$-expansion of $ l_{\beta} $ is not periodic with odd period. We recall in the following lines, the construction of this code. If for any integer $ n \in \N^{\times} $, $ d_{2n} < d_1 $, then $S_{\beta}$ is coded by the language of words ending by a string of the type $ d_1 \cdots d_k j $ with $(-1)^{k+1}(d_{k+1} - j)< 0$. More generally, we define $(d_i)_{i \geq 1} $ thanks to two sequences $ (n_i)_{i \geq 1} $ and $(p_i)_{i \geq 1}$ such that for any integer $ i $, $ d_{2n_i-1+m} = d_m $ if $ 1 \leq m \leq p_i $ and $ (-1)^{2n_i+p_i}(d_{2n_i+p_i}-d_{p_i+1})< 0$.
\begin{equation}
 d_1d_2\cdots = d_1 d_2 \cdots d_{2n_1-1}d_1 \cdots d_{p_1} d_{2n_1 + p_1} \cdots d_{2n_2-1}d_1 \cdots d_{p_2} d_{2n_2 + p_2} \cdots.
 \label{6}
\end{equation}
We assume that if $ n = 0$, $ d_1 \cdots d_n $ is the empty word (in this case, $d_1 \cdots d_n x = x $). We set:

\begin{equation*}
\Gamma_0 = \{ d_1\cdots d_nj = y \in L_{\beta} | \text{ $ n\in \N $, $ d_1 \cdots d_{n+1-i} \prec \sigma^{i}(y)$, $ \forall i$ , $ 0 \leq i \leq n $}\};
\end{equation*}
\begin{equation*}
D= \{  d_1 \cdots d_{2k+1}, k \in \N \};
\end{equation*}
\begin{equation}
\Delta_{0}^0 = \{ d_1 \cdots d_{2k-1} \vert 2n_i+p_i \leq 2k-1< 2n_{i+1}-1 \text{ and } i \in \N\}\label{Dood0}
\end{equation}
with $ n_0 = p_0 = 0 $; 
\begin{equation*}
E=\{ d_1\cdots d_{2n_i-1}; \text{ $ i \in \N^{\times}$ }\};
\end{equation*}
\begin{equation}
\C=\{xy \in L_{\beta} | x \in \{\varepsilon\}\cup \left(D^*\cap L_{\beta}\right),\hspace{2mm} y \in \Gamma_0 \} \label{C}; 
\end{equation}
\begin{equation}
\Delta_0 = \{ x y \in L_{\beta} | x \in \{\varepsilon\}\cup \left(E^*\cap L_{\beta}\right), \hspace{2mm} y \in \Delta_0^0 \}\label{Dood}.
\end{equation}
In what follows, we denote by $ \phi $ the morphism defined from $ \{0, 1\} $ to $ \{0, 1 \}^*$ by : 
\begin{equation}
\begin{aligned}
 \phi :  \{0 &, 1 \} &&\longrightarrow \{0, 1\}^* \\
       &0&&\longrightarrow \phi(0)=1 \\
       &1&&\longrightarrow \phi(1)=100.
\end{aligned}
\end{equation}
We set $ u_n = \phi^n(1) $, $ v_n = u_{n-1}u_{n-1} = \phi^n(00)$ and $ \gamma_n $ denotes the real number such that $ \frac{1}{\gamma_n} $ is the smallest real satisfying:
\begin{equation}
1=\frac{1}{\gamma_n^{l(u_n)}}+\frac{1}{\gamma_n^{l(v_n)}}.
\end{equation}
In fact, $\gamma_n$ is the algebraic integer of the polynomial $ X^{l_n}-X-1 $ where $ l_n = max(l(u_n),l(v_n))$ (see the proof of Proposition 5 in \cite{NguemaNdong20161}). Note that $ \gamma_0 $ is the golden ratio $ \frac{1+\sqrt{5}}{2} $. 
We have the following results: 

\begin{theo}\label{thIE}
 Let $ \beta < - 1$. Then $ (S_{\beta}, \sigma )$ is an intrinsically ergodic dynamical system. The maximal entropy measure is the Champernowne measure of a prefix recurrent positive code. 
 \end{theo}

\begin{theo}\label{pp1}
Let $ \beta < - 1 $. The $\beta$-shift endowed with the shift is coded if and only if $ \beta \leq -\frac{1+\sqrt{5}}{2}$ and the $ \beta $-expansion of the left endpoint $ l_\beta $ is not periodic with odd period. In all cases, the support of a measure of maximal entropy is coded by a recurrent positive (prefix or suffix) code. 
\end{theo}

\begin{theo}\label{th3}
Let $ \beta$ be a real number less than $ -1$, $d(l_{\beta}, \beta) = (d_i)_{i \geq 1} $ defined in \eqref{6} not periodic with odd period. Then, all measures with the maximal entropy have the same support.
\begin{itemize}
\item[(i)] If $ \beta< - \gamma_0$, then the support is $S_{\beta}$ which is coded by $ \C $;
\item[(ii)] If $ \beta =-\gamma_0 $, then the support is $S_{-\gamma_0}$ which is coded by $ \{1, 00\} $;
\item[(iii)] If $ -\gamma_0< \beta< - \gamma_1$, then the support is coded by the language $\Delta_0$;
\item[(iv)] If $-\gamma_n \leq \beta < -\gamma_{n+1} $, then there is $ x\in [-\gamma_0,-\gamma_1 [$ such that $ d(l_{\beta}, \beta) = \phi(l_x,x) $. The support of measures with the maximal entropy of $S_{\beta}$ is coded by $\Delta_n = \phi^n(\Delta_0)$, where $ \Delta_0 $ is a code of the support of  measures with the maximal entropy of $ S_x $.
\end{itemize} 
\end{theo}

\subsection{A positive recurrent code}

We start this part by proposing more details (given in \cite{NguemaNdong2019}) on the construction of $ \C $ and $ \Delta_n $, with $ n \in \N $.
We set $ B_i = d_1d_2\cdots d_{2n_i-1} $ and 
\begin{equation}
\Delta_{0}^1 = \{ B_{k_1}\cdots B_{k_m} X\vert p_{k_i}< 2n_{k_{i+1}}-1, X \in \Delta_{0}^0 \text{ and } l( X ) > p_{k_t} \} \label{Dood1},
\end{equation}
\begin{equation*}
\Delta_{0} = \begin{cases}
\Delta_{0}^0\cup \Delta_{0}^1 &\text{ if \eqref{6} is satisfied }\\
\{d_1\cdots d_{2k+1} | k \in \N \} & \text{ if $ d_{2i}< d_1 $, $ \forall i \in \N^{\times} $}
\end{cases},
\end{equation*}
\begin{equation*}
x \in \Gamma_0\Leftrightarrow \begin{cases}
                              x=d_1 \cdots d_{n}j, &\text{ with $n \in \N $}\\
															(-1)^{n+1}(d_{n+1}-j)< 0, &\text{ with $ 0\leq j<d_1 $},\\
															 2n_i+p_i < n \leq 2n_{i+1}-1, &\text{ $i\in \N$, 
															$n_0=p_0 = 0$}; \\
(-1)^{p_i}d_{p_i+1}>(-1)^{p_i}j>(-1)^{p_i}d_{2n_i+p_i} &\text{ if $ n=2n_i+p_i-1 $};
															\end{cases}
\end{equation*}
\begin{equation*}
x \in \Gamma_1 \Leftrightarrow \begin{cases}
                    x=B_{k_1}\cdots B_{k_m}y &\text{ with $y\in\Gamma_0$, $l(y\geq p_{k_m}+2$ };\\
								     k_1, \cdots, k_m \in \N^{\times}\\  
										p_{k_i}<2n_{k_{i+1}}-1,  &\text{ $ 1\leq i \leq m-1 $ };
													     \end{cases}
\end{equation*}
\begin{equation*}
x\in \Gamma_1^{'} \Leftrightarrow \begin{cases}
x=  B_{k_1}\cdots B_{k_{m-1}} y\\
y\in \Gamma_0^{'} &\text{ with $ l(y)=2n_{k_m}+p_{k_m}$},\\
 p_{k_i} < 2n_{k_{i+1}}-1 &\text{ and } 1\leq i \leq m-2.
\end{cases}
\end{equation*}

\begin{equation}
\Gamma = \begin{cases}
\Gamma_0 \cup \Gamma_1 \cup \Gamma_1^{'} &\text{ if \eqref{6} is satisfied } \\
\{d_1 \cdots d_{n}j | (-1)^{n+1}(d_{n+1}-j)<0, 0 \leq j < d_1, n\in \N \} & \text{ if $ d_{2i}< d_1 $, $ \forall i \in \N^{\times} $}.
\end{cases}
\label{G}
\end{equation}
\begin{equation*}
\C = \{ x y | x \in \Delta_{0}^{*}, y \in \Gamma, l(y)\geq 2 \} \cup \Gamma.
\end{equation*}
Moreover, 
\begin{equation}
J(0) = \{ t \vert p_t < 2n_1-1 \}, \label{J0}
\end{equation}
and for all $ i \in \N^{\times}$, 
\begin{equation}
J(i) = \{ t \vert 2n_i-1 \leq p_t < 2n_{i+1}-1 \} \label{Ji};
\end{equation}
and we denote by $ \Delta_{i} $ the set of words $ x $ such that
\begin{equation}
\begin{cases} 
                    x = B_{t_1}\cdots B_{t_m}, \\
                    p_{t_k} \leq 2n_{t_{k+1}}-1, \\
                     t_m \in J(i), \\
                      t_k \not\in J(i) \\
                    \text{for $ k\neq m $ }, p_{t_m}< 2n_{t_1}-1.
                    \end{cases} \label{Di}
\end{equation}

By construction, $ \C $ and $ \Delta_{n}$ are codes. The language $ \C $ is prefix but $ \Delta_n $ is certainly suffix. It suffices to use a permutation on words of $ \Delta_{n} $ to obtain a prefix code. For instance words of $ \Delta_{0} $ are of the form $ d_1 x $. The set of words of the type $ x d_1 $ (with $ d_1x \in \Delta_{0}$ ) is prefix. It generates a free monoid having the same language with that of $ \Delta_{0}^{*}$. 

\begin{exple}\label{ex1}
Let $ \beta<0 $ be the algebraic integer of the polynomial 
\begin{equation}
P(X)= X^{15}+3X^{14}-X^{12}+2X^{11}-X^{6}-X^{5}+2X^{4}-X^{3}-2X^{2}+2X+1.
\end{equation}
Then, $d(l_{\beta}, \beta) = \cdot 2012121201200 \overline{21} $, \\
$ \Delta_1 = \{201, \overline{2012121}^k 20121, \overline{2012121}^k2012121201200\overline{21}^p| p, k\geq 0 \} $, \\
 $ \Delta_2 = \{2012121 \} $, $ \Delta_0 = \{2\}$, $ \Gamma_0 = \{0, 1, 21, 200 \} $, \\$ \Gamma_1 = \{ x 200 \vert x \in \Delta_1^* \} $ and $ \C = \{ x y\vert x \in \Delta_0^*, y \in \Gamma, l(y) \geq 2 \} \cup \Gamma $.
\end{exple}

\begin{exple}
Let $ \beta < -1 $ be the algebraic integer of the polynomial
\begin{equation}
P(X) = X^{14}+2X^{13}-2X^{12}+X^{11}+X^{10}-X^{9}+X^{8}-X^{7}+X^{6}-2X^{5}+X^{4}+X^{3}-2X^{2}+1.
\end{equation}
Then, $ d(l_{\beta},\beta) = \cdot 2012121201200 \overline{1} $. \\
 $ \Delta_1 = \{ 201, \overline{2012121}^k 20121 | k \geq 0 \} $, $\Delta_2 = \{2012121 \} $;\\
 $ \Delta_0^0 = \{2, 2012121201200(11)^k| k \geq 0\} $\\ $ \Delta_0^1 = \{xy2012121201200(11)^k| k \geq 0, x \in \Delta_1^{*}, y \in \Delta_2^* \} $,\\ 
$ \Gamma_0 = \{0,1, 21, 200, 2012121201200(11)^k10| k \geq 0 \} $; \\
$ \Gamma_1 = \{ x200, xy 2012121201200(11)^k10| k \geq 0, x \in \Delta_1^*, y \in \Delta_2^* \} $
\end{exple}


 \begin{prop}\label{ppp1}
  Let $ \beta <- 1 $ and $ d(l_\beta, \beta) = (d_i)_{i \geq 1}$. We denote by $ D $ the $\beta$-shift subset of admissible concatenations of words of the type $ d_1 \cdots d_{2n-1} $ and $m$ an ergodic measure on $ S_{\beta} $ with the maximal entropy. If $m$ has support $D$, then: 
 \begin{equation}
 h_m(S_{\beta}) \leq \log \frac{1+\sqrt{5}}{2}.
 \label{(26)}
 \end{equation}
 \end{prop}


\begin{preu}
Let $F^{\times}$ be the set of words on $ \{0, 1, \cdots, d_1\} $ which can be decomposed into products of beginnings of $ (d_i)_{i \geq 1}$ of odd length and $ F = F^{\times} \cup \{\varepsilon\} $. Then, $ D \subset F $. But
\begin{equation}
F^{\times} = d_1F \cup  (d_1d_2d_3F) \cup (d_1d_2d_3d_4d_5F) \cup \cdots.
 \end{equation}
 This implies that the number $f_n$ of words with length $n$ of $F$ satisfying:
\begin{align*}
  f_n &= f_{n-1} + f_{n-3} + f_{n-5} + \cdots. \\
      &= f_{n-1} + f_{n-2}.
 \end{align*}
And thus, $ \frac{1}{n}\log f_n $ tends to $ \log \frac{1+\sqrt{5}}{2} $. Since $ D \subset F $, we obtain the result.
\begin{equation}
 h_m(S_{\beta}) \leq \log \frac{1+\sqrt{5}}{2}.
\end{equation}

\end{preu}

\begin{prop}\label{pp2}
Let $ \beta \in [ -\gamma_{n}, -\gamma_{n+1} )$. Then $ d(l_\beta, \beta) $ is the image under $ \phi^n$ of the $x$-expansion of $l_x=\frac{x}{1-x}$, for some $ x $ satisfying $ -\gamma_1 > x \geq -\gamma_0$.
\end{prop}

\begin{preu}
 Consider $ \beta \in [ -\gamma_{n}, -\gamma_{n+1} ) $. Then:
\begin{equation*}
d(l_{-\gamma_n},-\gamma_n) \preceq d(l_\beta,\beta) \prec d(l_{-\gamma_{n+1}}, -\gamma_{n+1}). 
\end{equation*}
So, there exists $ n_1 $ such that $ d(l_\beta, \beta) $ starts with $ u_n(u_{n-1})^{2n_1}u_n $. The word $ (u_{n-1})^{2n_1} $ is the longest concatenation of $ u_{n-1}$ which follows $ u_n $ in an admissible sequence. But all sequences of length $ l( u_n(u_{n-1})^{2n_1}u_n )=2l(u_n)+2n_1l(u_{n-1}) $ are greater than $ u_n(u_{n-1})^{2n_1}u_n $ (in the meaning of the alternating order) and after $ u_n $ one must have an even number of $ u_{n-1} $. Hence, there exists a bounded sequence $ (n_i)_{i \geq 1} $, $ n_i \leq n_1 $ such that  
 \begin{equation}
  d(l_{\beta}, \beta) = u_n(u_{n-1})^{2n_1}u_n(u_{n-1})^{2n_2}u_n (u_{n-1})^{2n_3} \cdots = \phi^n(1(0)^{2n_1}1(0)^{2n_2} \cdots).
 \end{equation}
 Moreover, $ d(l_{-\gamma_n}, -\gamma_n) = u_n(u_{n-1}u_{n-1})^{\infty} = \phi^n(1(0)^{\infty})$.
 Hence the result follows.
\end{preu}

\begin{lem}\label{lm1}
Let $ \beta $ be a real number such that  $ -1> \beta \geq -\gamma_0 $. Then  $ Card (\Delta_0 ) \geq 2$ if and only if $ \beta <-\gamma_1 $, where $ Card( \Delta_0 ) $ is the cardinality of $ \Delta_0 $.
\end{lem} 

\begin{preu}
\begin{itemize}
\item Suppose $ Card( \Delta_0 ) \geq 2 $. Note that when $ -\gamma_0$  is less than $ \beta $, the length of any string of zeros which appears in the $\beta$-expansion of $ l_{\beta}$ is even. Moreover, remark that if $ \Delta_0 $ contains a word of length 3 (that is 100), then 10000 is admissible and thus $ d(l_{\beta}, \beta)\prec d(l_{-\gamma_1}, -\gamma_1)$. In fact $ d(l_{-\gamma_1}, -\gamma_1) = 100\overline{11} $ and $ 10000 \prec 10011$.

Now suppose that $ 100 \notin \Delta_0 $, $ d(l_{\beta}, \beta) $ starts with $ 10011 $. For $ u \in \Delta_0 $ with length at least $ 5 $, if $ 10011 \notin\Delta_0$, $ 10011u \in \Delta_0 $ (by definition of $ \Delta_0$). The longest sequence of zeros is 00. So, $u$ ends with a sequence of the type $100(1)^t$ for some integer $ t $. Since $ 1\in \Delta_{0}$ (all concatenations of words of $ \Delta_0$ are admissible), it follows that $ 100(11)^{\infty}$ is admissible, that is $ d(l_{\beta}, \beta) \prec 100\overline{11} $. Thus, from Proposition 2 of \cite{NguemaNdong20161}, $ \beta <- \gamma_1 $.

\item Suppose $ -\gamma_1 > \beta \geq -\gamma_0 $.
\begin{equation}
1(0)^{\infty} \prec d(l_{\beta}, \beta) \prec 100(11)^{\infty}.
\end{equation}
Then, $ d(l_{\beta}, \beta) = 100(11)^{t_1}00(1)^{t_2}00(1)^{t_3} \cdots $. If $ t_1 = 0$, $100 \in \Delta_0$. If $t_1 \neq 0$, $ 100(11)^{t_1} \in \Delta_0$. Consequently, $Card(\Delta_0) \geq 2$.
\end{itemize}
\end{preu} 

\begin{lem}\label{lm2}
Let $ \beta $ be a real number such that $ -\gamma_1 > \beta \geq -\gamma_0$. Then, the topological entropy of $ (\Delta_0^*, \sigma) $ is larger than that of $ (\Delta_n^*, \sigma) $, with $ n \geq 1$.
\end{lem}

\begin{preu}
From Lemma \ref{lm1}, $ Card(\Delta_0) \geq 2 $, that is, $ \Delta_0 \neq \{ d_1 \} $. Remark that if $ d_1 \cdots d_{2n_i-1} \notin L_{\Delta_0^*}$, then for all $ t \geq i $, $ d_1 \cdots d_{2n_t-1} \notin L_{\Delta_0^*}$. If such an integer $ i $ is minimal, $ d(l_{\beta}, \beta)$ is an infinite concatenation of two consecutive words (with respect to the alternating order) $ U_0 = d_1 \cdots d_{2n_i-1}$ and $ V_0 = d_1 \cdots d_{2n_i-2}(d_{2n_i-1}-1)0 $ or $ V_0 = d_1\cdots d_{2n_i-3}(d_{2n_i-2}+1) $ (see the proof of Theorem 2 of \cite{NguemaNdong2019}). In an infinite admissible sequence, $ U_0 $ and $ V_0$ are followed by $ U_0 $ or $ V_0$. Thus, we can find in $\Delta_0 $ a word $ x \neq d_1 $ with length less than $ l(U_0)$ and if $ \Delta_{i_0} $ is the language which contains $ U_0 $, one has $ \Delta_{i_0}^* \subset \{U_0, V_0 \}^{*} $. 

Let $ \log\beta_1 $ be the entropy of $ \{U_0, V_0\}^* $ endowed with the shift. One has
\begin{equation}
\begin{aligned}
1 &= \frac{1}{\beta_1^{l(U_0)}}+\frac{1}{\beta_1^{l(V_0)}}\\
  &= \sum\limits_{n \geq 0} \frac{1}{\beta_1^{nl(U_0)+l(V_0)}}\label{U0}
\end{aligned}
\end{equation}
The entropy of $ \Delta_{i_0}^*$ is less than $ \log\beta_1$ since $ \Delta_{i_0}^* \subset \{U_0, V_0 \}^{*} $.

We have seen that there is a word $ x $ in $\Delta_0 $ such that $ x \neq d_1$ and $l(x)< l(U_0)$. Then $ \{d_1, x \}^* \subset \Delta_0^* $. Let $ \log \beta_2 $ be the entropy of $ \{d_1, x \}^*$. We have
\begin{equation}
1 = \frac{1}{\beta_2}+\frac{1}{\beta_2^{l(x)}},
\end{equation}
Let $ d_{\beta_2}(1) $ be the $\beta_2$-expansion of 1. Then, $ d_{\beta_2}(1) = 1(0)^{l(x)-1}1 $. So, because $l(x)< l(U_0)$,  $ (0)^{l(V_0)-1}\overline{1(0)^{l(U_0)-1}}$ is an infinite word of the $\beta_2$-shift. This implies that $ \sum\limits_{n \geq 0} \frac{1}{\beta_2^{nl(U_0)+l(V_0)}}<1 $. Consequently, $\beta_1<\beta_2$ since $ \beta_1 $ is the largest real satisfying \eqref{U0} and the map $ z\mapsto \sum\limits_{n \geq 0} \frac{1}{z^{nl(U_0)+l(V_0)}}$ on $ \R_{+}^{*}$ is decreasing. Hence the result follows.
\end{preu}

\begin{rem}\label{nlem1}
If $ -\gamma_{n+1}> \beta \geq -\gamma_n $, $\Delta_k = \{ u_k \} $, for all $ k< n$ (see the proof of Lemma 7 of \cite{NguemaNdong2019}).
\end{rem}

\begin{lem}\label{lmm3}
Let $ \beta < -1 $. Then:
\begin{itemize}
\item if $ \beta < -\gamma_0 $, then $ \sum\limits_{x \in \C}\frac{l(x)}{|\beta|^{l(x)}} <+\infty $;
\item if $ -\gamma_n \leq \beta < -\gamma_{n+1}$, then $ \sum\limits_{x \in \Delta_n}\frac{l(x)}{|\beta|^{l(x)}} < +\infty $.
\end{itemize}
\end{lem}

\begin{rem}
The coefficients of the expansion in the formal power series of $ \dfrac{1}{\prod\limits_{k \geq 0}(1-\sum\limits_{x \in \Delta_{k}}z^{l(x)})} $ count the admissible concatenations of words  of the type $ d_1 \cdots d_{2k+1} $, with $k \in \N$.
\end{rem}

\begin{preu 5}

If $ H_n $ denotes the number of words of length $ n $ in $ L_{\beta}$, it follows from Proposition 1 of \cite{NguemaNdong20161} that:
\begin{equation}
H_n = \sum\limits_{k=1}^{n}(-1)^k(d_{k-1}-d_k)H_{n-k}+1.
\end{equation}
Using Theorem 2 and Theorem 3 of \cite{NguemaNdong2019}, in the sense of formal power series, we have the following equation:
\begin{equation}
1-\sum\limits_{n \geq 1} (-1)^n(d_{n-1}-d_n)z^n = (1+z)(1-\sum\limits_{x\in \C}z^{l(x)})\prod\limits_{k \geq 0}(1-\sum\limits_{x \in \Delta_{k}}z^{l(x)}). \label{lab1}
\end{equation}
The left power series vanishes at $ -\frac{1}{\beta} = \frac{1}{|\beta|}$ which is its smallest root in modulus. That is $ \sum\limits_{n \geq 1}\frac{d_{n-1}-d_n}{\beta^n} = 1$ and so: 
$ \sum\limits_{x\in \C}\frac{1}{|\beta|^{l(x)}} = 1 $
or there exists $ n \in \N $ such that $  \sum\limits_{x\in \Delta_{n}}\frac{1}{|\beta|^{l(x)}}=1 $. 

From Proposition \ref{ppp1}:
\begin{itemize}
\item[1.] If $ \beta <- \gamma_0 $, $S_{\beta}$ is coded by $ \C $. Moreover the entropy of the system generated by the language $ \{ d_1\cdots d_{2n+1}| n \in \N\} $ is $\log \gamma_0 $. Thus the subsystem of admissible sequences which are concatenations of words of the type $ d_1\cdots d_{2n+1}$ has entropy less than $ \log \gamma_0 $. Then 
\begin{equation}
 \prod\limits_{k \geq 0}(1-\sum\limits_{x \in \Delta_{k}}\frac{1}{|\beta|^{l(x)}}) \neq 0 \text{ and } 1-\sum\limits_{x\in \C}\frac{1}{|\beta|^{l(x)}} = 0. 
\end{equation}
\item[2.] If $ -\gamma_n \leq \beta <- \gamma_{n+1} $, $\C=\{0\}$, for $ i < n$, $ \Delta_i = \{u_i\} $ and for $ i > n $, $ \Delta_i^* \subset L_{\Delta_n^*}$ (see the proof of Lemma 7 of \cite{NguemaNdong2019}). Thus, the entropy of $ \Delta_n^* $ is greater than that of the sub-system of  concatenations of words of the sets $ \Delta_i $ with $ i \geq n+1 $. The coefficients of the expansion in the formal power series of $ \dfrac{1}{\prod\limits_{k \geq n+1}(1-\sum\limits_{x \in \Delta_{k}}z^{l(x)})} $ count admissible concatenations of words of the sets $ \Delta_i $ with $ i \geq n+1 $. Then $ 1-\sum\limits_{x \in \Delta_{n}}\frac{1}{|\beta|^{l(x)}} = 0$ and  $ \prod\limits_{k \geq n+1}(1-\sum\limits_{x \in \Delta_{k}}\frac{1}{|\beta|^{l(x)}}) \neq 0 $.
\end{itemize}
We now have
\begin{equation}
\begin{aligned}
\sum\limits_{x \C}\frac{1}{|\beta|^{l(x)}} = 1 &\text{ if $ \beta <- \gamma_0 $}\\
\sum\limits_{x \in \Delta_n }\frac{1}{|\beta|^{l(x)}} = 1 &\text{ if $ -\gamma_n \leq \beta < -\gamma_{n+1} $ }
\end{aligned}\label{eqCD2}
\end{equation}
Using the derivatives of the formal powers series in \eqref{lab1} and the relation \eqref{eqCD2}, one has:
\begin{equation}
\sum\limits_{n \geq 1}n\dfrac{(d_{n-1}-d_n)}{\beta^n} = 
\begin{cases}
(1+\frac{1}{|\beta|}) \sum\limits_{x \in \C}\frac{l(x)}{|\beta|^{l(x)}} \prod\limits_{k \geq 0} \left(1-\sum\limits_{x \in \Delta_k}\frac{1}{|\beta|^{l(x)}}\right) &\text{ if $ \beta <- \gamma_0 $ } \\
(1-\frac{1}{\beta^2})\sum\limits_{x\in \Delta_n}\frac{l(x)}{|\beta|^{l(x)}}\prod\limits_{k \neq n} \left(1-\sum\limits_{x \in \Delta_k}\frac{1}{|\beta|^{l(x)}}\right) &\text{ if $ -\gamma_n \leq \beta < -\gamma_{n+1} $}
\end{cases} 
\end{equation}
Since $ (d_{n-1}-d_n)_{n \geq 1} $ is bounded, it follows that 
\begin{equation}
\begin{aligned}
\sum\limits_{x \C}\frac{l(x)}{|\beta|^{l(x)}} <+\infty &\text{ if $ \beta <- \gamma_0 $ } \\
\sum\limits_{x \in \Delta_n }\frac{l(x)}{|\beta|^{l(x)}} < +\infty &\text{ if $ -\gamma_n \leq \beta < -\gamma_{n+1} $ }.
\end{aligned}
\end{equation}
\end{preu 5}

\begin{preu}of Theorem \ref{pp1}
 \begin{itemize}
\item If $ \beta < -\gamma_0$, $ S_{\beta}$ (or the support of measures with maximal entropy) is coded by $ \C $ (see \cite{FloarXiv,NguemaNdong2019}). Then $ L_{\Delta_n^*}\subset L_{\C^*}$. The topological entropy of $ \C^* $ is larger than that of $ \Delta_n^*$, for all $ n $. In this case, it follows from \eqref{lab1} that $ \sum\limits_{x \in \C} \frac{1}{|\beta|^{l(x)}}= 1$.
\item If $ -\gamma_1 > \beta \geq -\gamma_0$, $ \C = \{ 0 \} $ and from Lemma \ref{lm2}, $ \sum\limits_{x \in \Delta_0} \frac{1}{|\beta|^{l(x)}}= 1$ The support of measures with maximal entropy is coded by $ \Delta_0 $.
\item If $ -\gamma_{n+1} > \beta \geq -\gamma_n $, it follows from Proposition \ref{pp2}, $ d(l_{\beta}, -\beta)$ is the image by $ \phi^n $ of a word taken between $ 1(0)^{\infty}$ and $ 100(11)^{\infty}$. Such a word is a concatenation of $1 $ and $ 00 $. In $[-\gamma_{n}, -\gamma_{n+1})$, the alphabet changes from $ \{1, 00\} $ to $ \{u_n, u_{n-1}u_{n-1} \}$. Since $ d_1 = 1 \in \Delta_0$, $ u_n = \phi^n(1)$ belongs to the language which codes the support of measures with maximal entropy. The language $ \Delta_i $ which contains $ u_n $ is $ \Delta_n$. It is the code of the support of measures with maximal entropy. Thus
$ \sum\limits_{x\in \Delta_n}\frac{1}{|\beta|^{l(x)}}=1$.
\end{itemize}
\end{preu}

\begin{preu 3}

Theorem \ref{th3} is due to the proof of Theorem \ref{pp1}, Lemma \ref{lm2} and Proposition \ref{pp2}. 
\end{preu 3}

From now on, we denote by $ P $ the code of the support of measures with maximal entropy.  
\begin{equation}
P = \begin{cases}
    \C &\text{ if $ \beta < - \gamma_0 $ } \\
		\Delta_n &\text{ if $ -\gamma_{n}\leq \beta < -\gamma_{n+1}$}.
		\end{cases}
\end{equation}

\subsection{Gaps on the negative beta-shift}

The phenomenon of gaps on $ I_{\beta} = [l_{\beta}, r_{\beta} ) $ was closely studied in \cite{MR2974214}. In this section, we are going to have the same study on the $\beta $-shift.

\begin{defi}
A word $ v \in L_{\beta} $ is intransitive if there exists $ u \in L_\beta$ such that for any $ w $ in $ L_{\beta}$, $ uvw \not\in L_{\beta} $. 
\end{defi}
We can see an intransitive word as a word which does not belong to the language of the support of a measure with maximal entropy. 

The following result is obvious.
\begin{prop}\label{pp5}
Let $ \beta $ be a real number such that
\begin{equation}
  d(l_{\beta}, -\beta) = u_n(u_{n-1})^{2k_1}u_n(u_{n-1})^{2k_2}u_n (u_{n-1})^{2k_3} \cdots.
 \end{equation}
An admissible word is intransitive if it contains one of the following sequences:
\begin{align*}
 &\sigma^i(u_{m-2}) u_{m-1}u_m &&\text{ with $ m > 0$, $ 0 \leq i < |u_{m-2} | $,}\\
 &\sigma^i(u_{m-1}) u_{m-1}u_{m-1}u_{m-1} &&\text{ with $ m \geq 0 $, $ 0 \leq i< |u_{m-1}|$, }\\
 &\sigma^i(u_{m-1}) u_{m-1} \cdots u_{n-2} u_{n-2} (u_{n-1})^{2k_1+1} u_n &&\text{ with $ m \geq 0 $, $ 0\leq i < |u_{m-1}|$.}
 \end{align*}
with $ u_{-1} = 0 $.
\end{prop}
The words listed in Proposition \ref{pp5} are forbidden in the language of the support of a measure with maximal entropy.

It is easy to see that an admissible word $ x $ starting with $ \sigma^i(u_k) $ contains an intransitive word if and only if it is taken between
\begin{equation}
\sigma^i(u_k)u_k u_{k+1}u_{k+1}\cdots u_{n-2} u_{n-2} (u_{n-1})^{2k_1}u_n (u_{n-1})^{2k_2}\cdots = \sigma^i(d(l_{\beta}, -\beta))
\end{equation}
and\begin{equation}
\sigma^i(u_k) u_{k+1}u_{k+1}\cdots u_{n-2} u_{n-2} (u_{n-1})^{2k_1}u_n (u_{n-1})^{2k_2}\cdots = \sigma^{|u_k|+i}(d(l_{\beta}, -\beta)).
\end{equation}

\begin{theo}\label{InW1}
Let $ \mu $ be an ergodic measure on the symbolic system $ (X, T)$ and $ L $ its language.
Consider two words $ u $ and $ t $ of  $ L $ such that $ \forall a \in L $, $ uat \not\in L$ ( $ t$ is intransitive if there is such a word $ u $). Then, $ \mu(_0[t]) = 0 $ or $ \mu(_0[u]) = 0$.
\end{theo}

\begin{preu}
When a measure $ \mu $ is ergodic, almost every point is generic (see Proposition (5.9) of \cite{MR0457675}). Thus, if $\mu(_0[u])$ and $ \mu(_0[t]) $ are not equal to zero and if $(x_n)_{n\in\Z}$ is generic for $ \mu $, there exist infinitely many words $ u $ and $ t $ in the sequence $ (x_n)_{n \geq 1 }$ and thus a word $ a $ in $ L $ such that $ u a t \in L $ (and a word $ b $ of $ L $ such that $ tbu \in L $).
\end{preu}

We deduce the following theorem.
\begin{theo}\label{InW2}
Let $\beta<-1 $ and $ \mu $ be a measure with maximal entropy on the negative $ \beta$-shift; $ \mu (_0[x])>0 $ whenever $ x $ can be decomposed into product of words of the recurrent positive prefix (or suffix) code of its support. Let $ t $ be an intransitive word of $ L_{\beta}$. Then $\mu(_0[t]) = 0$. 
\end{theo}

\subsection{Intrinsic ergodicity}

The existence of an ergodic measure on the tower implies that the set of ergodic measures on $S_{\beta}$ is non-empty. Indeed, there exists an invariant probability $\overline{\mu}$ with entropy $\log |\beta|$ on the tower associated to the code $\C$ inducing the Bernoulli probability $ \mu$ on the set of infinite words of the free monoid generated by the code and defined by:
\begin{equation}
 \mu ( [ x ] ) = \frac{1}{|\beta|^{l( x )}}, \text{ with $ x \in \C $}.
\end{equation}
Given a prefix (or suffix) code $P$, we denote by $W(P)$ the set of infinite sequences which can be decomposed into a product of words of the code $ P $. Let $\nu$ be the map from $ \F(S_{\beta})$ (subsets of $ S_{\beta}$) into $ [0, 1] $ which coincides with $\mu$ on all subsets of $ W(P)$ and zero on all subsets of the complement of $W(P)$ in $S_{\beta}$.

\begin{equation}
 \nu (B) = \begin{cases}
        \mu (B) & \text{ if $ B \subset W(P) $ }  \\
        0       & \text{ if $ B \subset S_{\beta} / W(P)$ } 
       \end{cases}
\end{equation}
$ \nu $ is an ergodic probability on $S_{\beta} $.

The existence of a recurrent positive (prefix or suffix) code implies the uniqueness of the measure with entropy $\log |\beta|$ on the tower associated to this prefix (or suffix code) code. However, there is one fact that needs to be taken account: all words cannot be written as concatenation of wors of the code. For instance, when  $\beta$ is less than $ -\frac{1+\sqrt{5}}{2} $, an infinite conctenation of the beginings with oddd length of the $\beta$-expansion of $\frac{\beta}{1-\beta} $ does not belong to $W(P)$ (for example $ \overline{ d_1 } $, $ d_1 \overline{d_1d_2d_3}, \cdots $). Also, the words ending with $ d(l_\beta, \beta) $ cannot be decomposed in $\C$ (we have for example the words $ d_1 d_1d_2d_3 \cdots$ , $ d_1d_2d_3 d(l_\beta, -\beta), \cdots $). Moreover, $ \cdot x_1 x_2\cdots $ and $ 0 \cdot x_1 x_2 \cdots $ denote the same $\beta$-representation in the $\beta$-shift. 

\begin{itemize}
\item[•] Suppose that $ \beta < -\dfrac{1+\sqrt{5}}{2} $ and $ d(l_{\beta}, \beta) $ is not periodic with odd period. 

Let  $ G $ be the set of sequences $ (x_i)_{i \geq n } $ (with $ n \leq 0 $) such that $ x_0x_1\cdots $ is on the form $X d(l_{\beta}, \beta)$ where $ X $ is the empty word or an admissible concatenation of words of $ \Delta_i $. 

We denotes by $ \tilde{\Delta_i^{\N}} $ the set of sequences $ (y_i)_{i \geq n } $ (with $ n \leq 0 $) such that $ y_0y_1 \cdots $ is an admissible concatenation of words of $ \Delta_i $.
\begin{equation}
S_{\beta} = \left( \underset{x \in \C}{\bigcup}[x]\right) \cup \left( \underset{i \geq 0}{\bigcup} \tilde{\Delta_i^{\N}} \right)  \cup G \label{s1}.
\end{equation}
 Since, the support of all measures with maximal entropy is coded by $ \C $, if $ \mu $ is one of these measures, one has
\begin{equation*}
 \mu \left( \left( \underset{i \geq 0}{\bigcup} \tilde{\Delta_i^{\N}} \right)  \cup G \right)  = 0 \text{ and }
\mu(S_{\beta}) = \mu \left( \underset{x \in \C}{\bigcup}[x]\right).
\end{equation*}

\item[•] If $ -\gamma_n \leq \beta < -\gamma_{n+1} $ and $ d(l_{\beta}, \beta) $ is not periodic with odd period, we denote by $ I_n $ the set of intransitive words. Then,
\begin{equation}
S_{\beta} = \left( \underset{ x \in I_n }{\bigcup}[x] \right) \cup \left( \underset{x \in \Delta_n}{\bigcup}[x] \right) \cup\left( \underset{i\geq n+1}{\bigcup} \tilde{\Delta^{\N}_i} \right) \cup G\label{s2}. 
\end{equation}
Since all measures with maximal entropy have the same support coded by $ \Delta_n $, if $ \mu $ is one of these measures, one has  
\begin{equation*}
\mu\left( \left( \underset{ x \in I_n }{\bigcup}[x] \right) \cup\left( \underset{i\geq n+1}{\bigcup} \tilde{\Delta^{\N}_i} \right) \cup G \right) = 0 \text{ and } \mu \left( \left( \underset{x \in \Delta_n}{\bigcup}[x] \right)\right) = \mu(S_{\beta}).
\end{equation*}
When $d(l_{\beta}, \beta) $ is periodic with odd period, the study can be done on $\tilde{S}_{\beta}$ since both systems $ S_{\beta}$ and $\tilde{S}_{\beta}$ have the same entropy. And thus $ S_{\beta}\setminus \tilde{S}_{\beta} $ is negligible with respect to  measure with maximal entropy.

\item[•] For $ d(l_{\beta}, \beta) $ periodic with odd period, we have
\begin{equation}
S_{\beta} = \tilde{S}_{\beta} \cup T \label{s3}
\end{equation}
where $ T $ is the set of admissible sequences ending with $ d(l_{\beta}, \beta) $.
\end{itemize}

From the previous arguments, we have seen that, for the study of intrinsic ergodicity, it suffices to concentrate our attention on cylinders carried by words of $ P $ which codes the support of measures with maximal entropy. And we can neglect all other sub-sets of $ S_{\beta}$.

We define on the tower $ (\Omega, T) $ of $ P $ the map $ f $ by:
\begin{equation}
 f( (x_n)_{n \in \Z}, i) =  (y_n)_{n \in \Z} 
 \label{(2)}
\end{equation}
where $ x_k \in \C $, $ y_i \in \{0, 1, \cdots, d_1 \} $, $ y_0 $ denotes the $i$-th letter of $x_0$. Let $ x $ such that  
\begin{equation}
x = \cdots x_{-m}\cdots x_{-1}x_0. x_1 x_2 \cdots x_m \cdots = \cdots z_{-n}\cdots z_{-1} z_0 z_1 \cdots z_n \cdots.
\end{equation}
In fact, $ \cdots x_{-n} \cdots x_{-1} x_0 \cdots x_n \cdots $ and $ \cdots z_{-m}\cdots z_{-1}z_0. z_1 \cdots z_m \cdots $ are two writings of $ x $. In first case, $ x $ is a word of the free monoid generated by $P$. In the second case, $ x $ is viewed as a word of $ \A^*$. Thus
\begin{equation}
f((x_n)_{n \in \Z}, i ) = \sigma^{i}((z_n)_{n \in \Z}) = (z_{n+i})_{n \in \Z}.
\end{equation}
The map $ f $ is one to one. Moreover, it is easy to see that $ f\circ T = \sigma \circ f $.
 
Now, we have all ingredients needed for proving Theorem \ref{thIE}.

 \begin{preu}of Theorem\ref{thIE}
 
 All measures with maximal entropy have the same support. 
Since $ f $ is one to one and $ f\circ T = \sigma \circ f $, each $ \sigma$-invariant measure $ \mu $ on $ W(P) $ generates a measure $ \mu \circ f $  on $ \Omega $ with the same entropy. We have seen that the code $ P $ is recurrent positive. Then, there is a unique measure with entropy $ \log |\beta| $ on $\Omega$. This implies the existence of a unique measure with entropy $ \log |\beta|$ on $ W(P) $.

The restrictions on $W(P)$ of measures with maximal entropy on $ S_{\beta}$ have entropy $\log |\beta|$. Then, they coincide on $ W(P) $. Therefore, the measure with maximal entropy on $ S_{\beta}$ is unique.
\end{preu}

After proving the intrinsic ergodicity, let us determine the measure with maximal entropy (denoted by $\mu_{\beta}$) on cylinders carried by words of the code $ P $. 
 
Any invariant probability $ \nu $ on $(P^{\Z}, \sigma_P)$ with finite average length $l(P, \nu) $ is induced by a unique invariant probability measure $ \overline{\mu} $ of $(\Omega, T ) $ (see \cite{MR858689}). The link between the entropies of the two measures is given by the Abramov formula:
\begin{equation}
h(\overline{\mu})l(\nu, P) = h(\nu).\label{abra}
\end{equation} 
We know that $P^{\Z}$ is identified with the base $P^{\Z} \times \{1\}$ of $\Omega$ (see \cite{MR858689}) and for a Borel subset $B$ of $ P^{\Z} $, 
\begin{equation}
 \nu (B) = \dfrac{ \overline{\mu}( B \times \{1\} )}{\overline{\mu}(P^{\Z} \times \{1\} )}.\label{indm}
\end{equation}
Since $ P $ is recurrent positive, there is a unique measure $\overline{\mu}$ with entropy $ \log|\beta| $ on the tower $(\Omega, T) $ which induces the unique invariant probability measure $\nu $ on $P^{\Z}$ such that:
\begin{equation}
\nu([x])= \frac{1}{|\beta|^{l(x) }} \text{ where $ x \in P $ }.
\end{equation}

So, for $ x \in P $, 
\begin{equation}
 \nu ([x] ) \overline{ \mu } (P^{\Z}\times \{1\}) = \overline{ \mu }( [ x ] \times \{i\} ).
 \label{(27)}
\end{equation}
Since $ \Omega = \bigcup\limits_{ x \in P } \left(\underset{i=1}{\overset{l( x )}{\cup}} [ x ] \times \{ i \} \right) $, one has
\begin{equation}
\begin{aligned}
 1 &= \sum \limits_{ x \in P} \underset{i = 1}{\overset{ l( x ) }{\sum}}\overline{\mu}\left([x] \times \{i\} \right) \\
   &= \sum \limits_{ x \in P } \underset{i = 1}{\overset{ l( x ) }{\sum}}\overline{\mu}\left(T^{-i+1}( [x]\times\{i\}) \right) \\
   &= \sum \limits_{ x \in P } \underset{i = 1}{\overset{ l( x ) }{\sum}}\overline{\mu}\left([x] \times \{1\} \right) \\
   &= \sum \limits_{ x \in P} l( x ) \overline{\mu}\left([x] \times \{1\} \right).
\label{(28)}
\end{aligned}
\end{equation}
Moreover, the average length of $P$ with respect to the measure $\nu$ is:
 \begin{equation}
\begin{aligned}
 l(P, \nu ) &= \sum \limits_{ x \in P } l( x ) \nu ([x]) \\
                   &= \dfrac{ 1 }{ \overline{ \mu }(P^{*}\times \{1\} ) } \sum \limits_{ x \in P} \overline{\mu}([x] \times \{1\} ) \\
                   &= \dfrac{ 1 }{ \overline{ \mu }(P^{*} \times \{1\} ) } \\
                   &= \sum \limits_{ x \in P } \dfrac{l(x) } {|\beta|^{l(x)}}.
\label{(29)}                   
\end{aligned}
\end{equation}
And thus, one has:
\begin{equation}
 \mu_{\beta} ( [ x ] ) = \left( |\beta|^{l( x )} \sum \limits_{ x \in P } \dfrac{ l(x)  } {|\beta|^{l(x)}} \right)^{-1}.
 \label{(30)}
\end{equation}

 \begin{theo}\label{Mix}
 The measure with maximal entropy of the negative beta-shift is mixing.
 \end{theo}
Before proving Theorem \ref{Mix}, let us show the following result:

\begin{prop}\label{GCD}
The G.C.D of lengths of words of codes previously constructed is 1.
\end{prop}
\begin{preu}
 For $ \beta \leq -\gamma_0$, the $\beta$-shift is coded (by $ \C $ if the inequality is strict and by $ \{1, 00\} $ if $ \beta = \gamma_0$). And also, the code contains at least one word of length 1. 
 
 If $ \beta \in [ -\gamma_0, -\gamma_1 )$, the support is coded by $\Delta_{0}$ which contains $ 1 = d_1 $.
 
 Therefore, consider $ \beta $ such that $ -\gamma_{n+1} > \beta \geq -\gamma_n $ with $ n > 1$. In this case, the support of the maximal entropy measure is coded by $ \Delta_{n}$. The words of this set are of the form 
 \begin{equation}
 u_nv_n^{n_1+1}u_{n}v_n^{n_2}\cdots u_n v_n^{n_{2k}}u_nv_n^t,
 \end{equation}
 with $ 0 \leq n_{2k+1}-1$ and $ 0 \leq k$. 
 
If $ n_1 \neq 0$, $ u_n $ and $ u_nv_n $ belong to $ \Delta_{n} $. The integer $ l( u_n ) $ and $ l( v_n ) $ are relatively prime since 
 \begin{equation}
 l( u_n ) = l( v_n ) + (-1)^n.
 \end{equation}
Thus, $ l( u_n ) $ and $ l( u_n v_n ) $ are relatively prime too.  

Note that $ v_n^{n_1+1}$ is the longest sequence of $ v_n $ in the support of a measure of maximal entropy. Thus, if $ n_1 = 0$, $n_2=0$, $ u_n$ and $ u_n v_n u_nu_n$ belong to the code. But $ l( u_n ) $ and $ l( u_n v_n u_nu_n )$ are relatively prime. It follows that, for all $ \beta <- 1$, the G.C.D of lengths of words belonging to the code of the support of the maximal entropy measure is 1. 
\end{preu}

An immediate consequence of the Proposition \ref{GCD} is that the restriction of the measure with maximal entropy on its support is mixing. Note that, if $ x $ is an intransitive word, $ [x] $ is $\sigma$-invariant. And then, for all $ n $ and $ y $ in the code of support,
\begin{equation}
\sigma^{-n}[x]\cap [y] = \varnothing 
\end{equation}
Thus
\begin{equation}
\lim\limits_{n\rightarrow +\infty} \mu(\sigma^{-n}[x]\cap [y]) = 0 = \mu([x])\mu([y])\label{eqmix}
\end{equation}
since $\mu([x]) = 0 $. Moreover, for all $ n $, $\sigma^{-n}[y]\cap [x] \subset [x ]$ and then
\begin{equation}
\lim\limits_{n\rightarrow +\infty} \mu(\sigma^{-n}[y]\cap [x]) = 0 = \mu([y])\mu([x]).
\end{equation}
If now, $ x $ and $ y $ are both intransitive words, $\sigma^{-n}[x] \cap [y]$ is negligible with respect to the measure with maximal entropy. Then \eqref{eqmix} is also satisfied. This proves Theorem \ref{Mix}.

\section*{}

In summary, we have seen that for each case studied, there exists a unique $\sigma$-invariant measure of maximal entropy. Considering the one side $\beta$-shift, the results remain valid. If $ d(l_{\beta}, \beta) $ is periodic with odd period, $ S_{\beta}$ and $\tilde{S}_{\beta}$ have the same entropy. The negative $\beta$-shift $ S_{\beta}$ is the union of $\tilde{S}_{\beta}$ which is intrinsically ergodic and the $\sigma$-invariant sub-set of words ending with $ d(l_{\beta}, \beta) $. When $ \beta $ is between $ -\frac{1+\sqrt{5}}{2} $ and -1, the system is not transitive. But $S_{\beta}$ remains intrinsically ergodic. In \cite{article}, an example of sub-system of $ S_{\beta}$ not intrinsically ergodic is given by: $ X = \{1^{\infty} \} \cup \{1^n2^{\infty} : n\geq 1 \} \cup \{2^{\infty} \}$. This sub-shift corresponds to $\{0^{\infty} \}\cup \{ 0^n1^{\infty}: n \geq 1 \} \cup \{1^{\infty} \} $ according to our definition of the negative $\beta$-transformation. It is easy to see that this sub-shift is contained in all negative $\beta$-shift. This example shows that in the intrinsically ergodic dynamical system, we can find sub-systems which do not have this property. But it is necessary to attach the condition to this sub-system to have an entropy strictly less than the entropy of the system.





\bibliographystyle{unsrt}
\bibliography{Mabiblio}

\end{document}